\documentclass[11pt]{article}
\usepackage{amsmath,amsthm,amssymb,amsfonts}
\usepackage{graphicx}

\newtheorem{theorem}{Theorem}[section]

\newtheorem{lemma}[theorem]{Lemma}
\newtheorem{proposition}[theorem]{Proposition}

\newtheorem{remark}[theorem]{Remark}

\begin{document}

\title{Multiple nonradial solutions for a nonlinear elliptic radial problem: an
improved result}
\author{{\small \smallskip }Sergio Rolando\smallskip \\
\textit{\small Dipartimento di Matematica e Applicazioni}\\
\textit{{\small Universit\`{a} degli Studi di Milano bicocca}}, \textit%
{\small Via Cozzi 53, 20125 Milano}, \textit{\small Italy}\\
\textit{\small e-mail: }{\small sergio.rolando@unito.it}}
\date{}
\maketitle

\begin{abstract}
We obtain an improved version of a recent result concerning the existence of
nonnegative nonradial solutions $u\in D^{1,2}(\mathbb{R}^{N})\cap L^{2}(%
\mathbb{R}^{N},\left| x\right| ^{-\alpha }dx)$ to the equation 
\[
-\triangle u+\displaystyle\frac{A}{\left| x\right| ^{\alpha }}u=f\left(
u\right) \quad \text{in }\mathbb{R}^{N},\quad N\geq 3,\quad A,\alpha >0, 
\]
where $f$ is a continuous nonlinearity satisfying $f\left( 0\right) =0$%
.\bigskip

\smallskip

\end{abstract}

\section{Introduction and main result}

In this paper we consider the following semilinear elliptic problem: 
\begin{equation}
\left\{ 
\begin{array}{l}
\begin{array}{ll}
\medskip -\triangle u+\displaystyle\frac{A}{\left| x\right| ^{\alpha }}%
u=f\left( u\right) ~ & \text{\textrm{in }}\mathbb{R}^{N},~N\geq 3
\end{array}
\\ 
\begin{array}{l}
u\in H_{\alpha }^{1}\setminus \left\{ 0\right\} ,\quad u\geq 0
\end{array}
\end{array}
\right.   \tag*{$\left({\mathcal P}\right) $}
\end{equation}
where $A,\alpha >0$ are real constants, $f:\mathbb{R}\rightarrow \mathbb{R}$
is continuous and satisfies $0<f\left( s\right) \leq \left( \mathrm{const.}%
\right) s^{p-1}$ for some $p>2$ and all $s>0$ (hence $f\left( 0\right) =0$),
and $H_{\alpha }^{1}$ is the natural energy space related to the equation,
i.e., 
\[
H_{\alpha }^{1}:=\left\{ u\in D^{1,2}(\mathbb{R}^{N}):\int_{\mathbb{R}^{N}}%
\frac{u^{2}}{\left| x\right| ^{\alpha }}dx<\infty \right\} .
\]
Here and in the rest of the paper, $D^{1,2}(\mathbb{R}^{N})$ is the usual
Sobolev space, which identifies with the completion of $C_{\mathrm{c}%
}^{\infty }(\mathbb{R}^{N})$ with respect to the $L^{2}$ norm of the
gradient.

We are interested in \emph{weak solutions} to $\left( \mathcal{P}\right) $,
i.e., functions $u\in H_{\alpha }^{1}\setminus \left\{ 0\right\} $ such that 
$u\geq 0$ almost everywhere in $\mathbb{R}^{N}$ and 
\begin{equation}
\int_{\mathbb{R}^{N}}\nabla u\cdot \nabla v\,dx+\int_{\mathbb{R}^{N}}%
\displaystyle\frac{A}{\left| x\right| ^{\alpha }}uv\,dx=\int_{\mathbb{R}%
^{N}}f\left( u\right) v\,dx\quad \text{for all }v\in H_{\alpha }^{1}.
\label{weak sol}
\end{equation}

As is well known, problems like $\left( \mathcal{P}\right) $ arise in many
branches of mathematical physics, such as population dynamics, nonlinear
optics, plasma physics, condensed matter physics and cosmology (see e.g. 
\cite{Fife,BFmonograph,YangY}). In this context, $\left( \mathcal{P}\right) $
is a prototype for problems exhibiting radial potentials which are singular
at the origin and/or vanishing at infinity (sometimes called the \emph{zero
mass case}; see e.g. \cite{Meder, BPR}).

\begin{figure}
\centering
\includegraphics{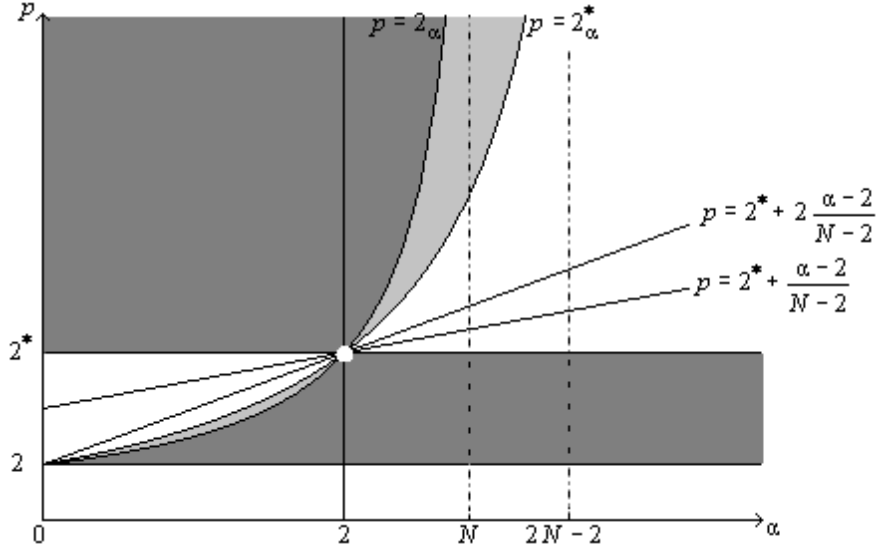} 
       \caption{Regions of nonexistence of solutions (light gray), and existence (white with $p>2$) and nonexistence (dark gray) of radial solutions.} 
								\label{pic}
\end{figure}

Although it can be considered as a quite recent investigation, the study of
problem $\left( \mathcal{P}\right) $ has already some history and,
currently, the problem of existence and nonexistence of \emph{radial}
solutions has been essentially solved, through various subsequent
contributions, in the \emph{pure-power case} $f\left( u\right) =\left|
u\right| ^{p-2}u$, where the results obtained rest upon compatibility
conditions between $\alpha $ and $p$. The first results in this direction
are probably due to Terracini \cite{Terracini}, who proved that $\left( 
\mathcal{P}\right) $ has no solution if 
\[
\left\{ 
\begin{array}{l}
\alpha =2 \\ 
p\neq 2^{*}
\end{array}
\right. \quad \text{or}\quad \left\{ 
\begin{array}{l}
\alpha \neq 2 \\ 
p=2^{*}
\end{array}
\right. ,\qquad 2^{*}:=\frac{2N}{N-2},
\]
and explicitly found all the radial solutions of $\left( \mathcal{P}\right) $
for $\left( \alpha ,p\right) =\left( 2,2^{*}\right) $. As usual, $2^{*}$
denotes the critical exponent for the Sobolev embedding in dimension $N\geq 3
$. The problem was then considered again in \cite{Co-Cr-Par}, where the
authors proved that $\left( \mathcal{P}\right) $ has at least a radial
solution if 
\[
\left\{ 
\begin{array}{l}
\smallskip 0<\alpha <2 \\ 
2^{*}+\frac{\alpha -2}{N-2}<p<2^{*}
\end{array}
\right. \quad \text{or}\quad \left\{ 
\begin{array}{l}
\smallskip \alpha >2 \\ 
2^{*}<p<2^{*}+\frac{\alpha -2}{N-2}
\end{array}
\right. ,
\]
while it has no solution if 
\[
\left\{ 
\begin{array}{l}
0<\alpha <2 \\ 
p>2^{*}
\end{array}
\right. \quad \text{or}\quad \left\{ 
\begin{array}{l}
\alpha >2 \\ 
2<p<2^{*}
\end{array}
\right. .
\]
The existence and nonexistence results of \cite{Co-Cr-Par} were subsequently
extended in \cite{BRpow}, by showing that $\left( \mathcal{P}\right) $ has
no solution also if 
\[
\left\{ 
\begin{array}{l}
0<\alpha <2 \\ 
2<p\leq 2_{\alpha }
\end{array}
\right. \quad \text{or}\quad \left\{ 
\begin{array}{l}
2<\alpha <N \\ 
p\geq 2_{\alpha }
\end{array}
\right. ,\qquad 2_{\alpha }:=\frac{2N}{N-\alpha },
\]
and obtaining at least a radial solution for every pair $\left( \alpha
,p\right) $ such that 
\[
\left\{ 
\begin{array}{l}
\smallskip 0<\alpha <2 \\ 
2^{*}+2\frac{\alpha -2}{N-2}<p<2^{*}
\end{array}
\right. \quad \text{or}\quad \left\{ 
\begin{array}{l}
\smallskip \alpha >2 \\ 
2^{*}<p<2^{*}+2\frac{\alpha -2}{N-2}
\end{array}
\right. .
\]
A further extension of this existence condition was found in \cite
{Su-Wang-Will-2,Su-Wang-Will-p}, where the authors proved that $\left( 
\mathcal{P}\right) $ has at least a radial solution for all the pairs $%
\left( \alpha ,p\right) $ satisfying 
\begin{equation}
\left\{ 
\begin{array}{l}
0<\alpha <2 \\ 
2_{\alpha }^{*}<p<2^{*}
\end{array}
\right. \quad \text{or}\quad \left\{ 
\begin{array}{l}
2<\alpha <2N-2 \\ 
2^{*}<p<2_{\alpha }^{*}
\end{array}
\right. \quad \text{or}\quad \left\{ 
\begin{array}{l}
\alpha \geq 2N-2 \\ 
p>2^{*}
\end{array}
\right. ,\qquad 2_{\alpha }^{*}:=2\frac{2N-2+\alpha }{2N-2-\alpha }.
\label{sww}
\end{equation}
Finally, the problem of radial solutions in the left open cases was solved
in \cite{BGRnonex} and \cite{Catrina}, where, respectively, it was proved
that the problem has no radial solutions for both 
\[
\left\{ 
\begin{array}{l}
0<\alpha <2 \\ 
2_{\alpha }<p\leq 2_{\alpha }^{*}
\end{array}
\right. \quad \text{and}\quad \left\{ 
\begin{array}{l}
2<\alpha <2N-2 \\ 
2_{\alpha }^{*}\leq p<2_{\alpha }
\end{array}
\right. .
\]
All these results are portrayed in the picture of the $\alpha p$-plane given
in Fig.\ref{pic}, where nonexistence regions are shaded in gray
(nonexistence of radial solutions) and light gray (nonexistence of solutions
at all, which includes both the lines $p=2^{*}$ and $p=2_{\alpha }$ except
for the pair $\left( \alpha ,p\right) =\left( 2,2^{*}\right) $), whereas
white color (of course above the line $p=2$) means existence of radial
solutions.

All the above results on the existence of radial solutions are obtained by
variational tecniques which can be extended in a standard way (see \cite
{BRpow,Su-Wang-Will-p}) to general continuous nonlinearities $f:\mathbb{R}%
\rightarrow \mathbb{R}$ satisfying the so-called \emph{Ambrosetti-Rabinowitz
condition} (i.e., assumption $(\mathbf{h}_{2})$ below) and $0<f\left(
s\right) \leq \left( \mathrm{const.}\right) s^{p-1}$ for all $s>0$ and for
some $p>2$ such that the pair $\left( \alpha ,p\right) $ satisfies (\ref{sww}%
).

Concerning nonradial solutions, Terracini proved in \cite{Terracini} that
problem $\left( \mathcal{P}\right) $ with $N\geq 4$, $\alpha =2$ and $%
f\left( u\right) =\left| u\right| ^{2^{*}-2}u$ has at least a nonradial
solution for every $A$ large enough. This just concerns the point $\left(
\alpha ,p\right) =\left( 2,2^{*}\right) $ in Fig.\ref{pic}, which brings
Catrina to say, in the introduction of his paper \cite{Catrina}: ``\textit{%
Two questions still remain: whether one can find non-radial solutions in the
case when radial solutions do not exist, or in the case when radial
solutions exist}''.

A first result in this direction has been obtained in \cite{R}, for
nonlinearities satisfying the following assumption:

\begin{itemize}
\item[$(\mathbf{h}_{0})$]  there exist $p_{1},p_{2}>2$ and $M>0$ such that $%
0<f\left( s\right) \leq M\min \{s^{p_{1}-1},s^{p_{2}-1}\}$ for all $s>0.$
\end{itemize}

\noindent Assumption $(\mathbf{h}_{0})$ is the so-called \emph{double-power
growth condition} and seems to be typical in nonlinear problems with
potentials vanishing at infinity (see e.g. \cite
{Azz-Pomp,G-R,BGR-p2,BGR-p,BPR, Meder,BBR-2, R} and the references therein).
It obviously implies the single-power growth condition $0<f\left( s\right)
\leq Ms^{p-1}$ for all $p\in [p_{1},p_{2}]$ and $s>0$, but it is actually
more stringent than that whenever it is assumed with $p_{1}\neq p_{2}$.

\begin{figure}
\centering
\includegraphics{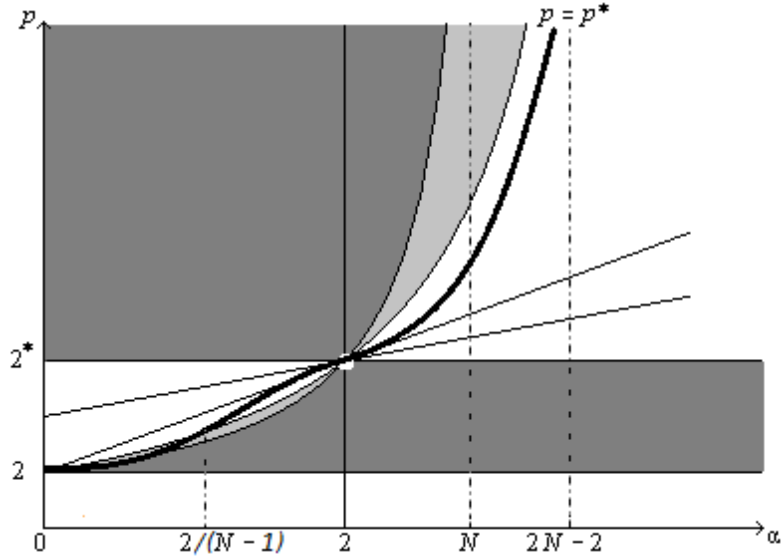} 
       \caption{The curve $p=p_{\alpha }^{*}$ in the same $\alpha p$-plane of Fig.\ref{pic}.} 
								\label{pic-star}
\end{figure}

The result of \cite{R} is the following. For every $\alpha \in \left(
0,2N-2\right) \setminus \left\{ 2\right\} $ and $p_{1},p_{2}>2$, define 
\[
p_{\alpha }^{*}:=\left\{ 
\begin{array}{lll}
2\frac{\alpha ^{2}(N-1)-2\alpha (N-1)+4N}{\alpha ^{2}(N-1)-2\alpha (N+1)+4N}%
\medskip & \text{if } & 0<\alpha <2 \\ 
2\frac{2N+2-\alpha }{2N-2-\alpha } & \text{if } & 2<\alpha <2N-2
\end{array}
\right. 
\]
and 
\begin{equation}
\nu :=\nu _{N,\alpha ,p_{1},p_{2}}:=\left\{ 
\begin{array}{lll}
\left\lceil 2\min \left\{ \frac{N-1}{\alpha },\frac{N-2}{2-\alpha }\frac{%
2^{*}-p_{1}}{p_{1}-2}\right\} -2N\left( \frac{1}{\alpha }-\frac{1}{2}\right)
\right\rceil -1\medskip & \text{if } & 0<\alpha <2 \\ 
\left\lceil 2\min \left\{ \frac{N-1}{\alpha },\frac{N-2}{\alpha -2}\frac{%
p_{2}-2^{*}}{p_{2}-2}\right\} \right\rceil -1 & \text{if } & 2<\alpha <2N-2
\end{array}
\right.  \label{n_alfa}
\end{equation}
where $\left\lceil \cdot \right\rceil $ denotes the ceiling function (i.e., $%
\left\lceil x\right\rceil :=\min \left\{ n\in \mathbb{Z}:n\geq x\right\} $).
Fig.\ref{pic-star} shows the curve $p=p_{\alpha }^{*}$ in the same $\alpha p$%
-plane of Fig.\ref{pic}. Note that $p_{\alpha }^{*}=2_{\alpha }^{*}$ if $\alpha =2/(N-1)$.

\begin{theorem}
\label{THM: main}Let $N\geq 4$ and $\alpha \in \left( 2/(N-1),2N-2\right) $, 
$\alpha \neq 2$. Let $f:\mathbb{R}\rightarrow \mathbb{R}$ be a continuous
function satisfying assumption $(\mathbf{h}_{0})$ with $2<p_{1}<p_{\alpha
}^{*}$ and $p_{2}>2^{*}$ if $\alpha \in \left( 2/(N-1),2\right) $, or $%
2<p_{1}<2^{*}$ and $p_{2}>p_{\alpha }^{*}$ if $\alpha \in \left(
2,2N-2\right) $. Assume furthermore that:

\begin{itemize}
\item[$(\mathbf{h}_{1})$]  the function $\text{$f$}$\text{$\left( s\right) /s
$} is strictly increasing on $\left( 0,+\infty \right) $

\item[$(\mathbf{h}_{2})$]  $\exists \mu >2$ such that the function $F$\text{$%
\left( s\right) /s^{\mu }$} is decreasing on $\left( 0,+\infty \right) $

\item[$(\mathbf{h}_{3})$]  $\exists \eta >2$ such that the function $F$\text{%
$\left( s\right) /s^{\eta }$} is increasing on $\left( 0,+\infty \right) $
\end{itemize}

\noindent where $F\left( s\right) :=\int_{0}^{s}f\left( t\right) dt$. Then
there exists $A_{*}>0$ such that for every $A>A_{*}$ problem $\left( 
\mathcal{P}\right) $ has both a radial solution and $\nu $ different
nonradial solutions.
\end{theorem}

Under the assumptions of Theorem \ref{THM: main}, $\nu $ turns out to be
strictly positive (see \cite[Lemma 5.2]{R}), so that at least one nonradial
solution actually exists. Unfortunately, the theorem does not encompass the
case of pure-power nonlinearities, since it always requires $p_{1}<p_{2}$ in
assumption $(\mathbf{h}_{0})$. Observe that assumption $(\mathbf{h}_{3})$ is
the well known \emph{Ambrosetti-Rabinowitz condition}, as the increasingness
of $F$\text{$\left( s\right) /s^{\eta }$} on $\left( 0,+\infty \right) $
amounts to $\eta F$\text{$\left( s\right) \leq f\left( s\right) s$} for all $%
s>0$.

The aim of this paper is to improve Theorem \ref{THM: main}, by proving the
following result.

\begin{theorem}
\label{THM: new}Let $N\geq 4$ and $\alpha \in \left( 2/(N-1),2N-2\right) $, $%
\alpha \neq 2$. Let $f:\mathbb{R}\rightarrow \mathbb{R}$ be a continuous
function satisfying $(\mathbf{h}_{0})$ with $p_{1},p_{2}$ as in Theorem \ref
{THM: main}. Assume furthermore that:

\begin{itemize}
\item[$(\mathbf{h}_{1}^{\prime })$]  the function $\text{$f$}$\text{$\left(
s\right) /s$} is increasing on $\left( 0,+\infty \right) $

\item[$(\mathbf{h}_{2}^{\prime })$]  $\exists \mu >2$ and $\exists s_{*}>0$
such that the function $F$\text{$\left( s\right) /s^{\mu }$} is decreasing
on $\left( 0,s_{*}\right) $.
\end{itemize}

\noindent Then the same conclusion of Theorem \ref{THM: main} holds true.
\end{theorem}

As for Theorem \ref{THM: main}, Theorem \ref{THM: new} does not concern
pure-power nonlinearities (owing to assumption $(\mathbf{h}_{0})$ with $%
p_{1}<p_{2}$), and we have $\nu \geq 1$ and $\lim
_{N\rightarrow\infty }\nu =+\infty $ (for $\alpha ,p_{1},p_{2}$ fixed).

Clearly Theorem \ref{THM: new} includes Theorem \ref{THM: main}, and in particular it removes the Ambrosetti-Rabinowitz condition from the assumptions.
Such inclusion is actually strict, for instance because Theorem \ref{THM:
new} applies to the nonlinearities 
\[
f\left( s\right) =\min \{\left| s\right| ,\left| s\right| ^{p_{2}-1}\}\quad 
\text{and}\quad f\left( s\right) =\frac{\left| s\right| ^{p_{2}-1}}{1+\left|
s\right| ^{p_{2}-2}}
\]
for any $p_{2}$ (and $p_{1}$) as in the theorem, while Theorem \ref{THM:
main} does not.

Theorem \ref{THM: new} will be proved in Section \ref{SEC: pf}. Our proof is
variational, since the weak solutions to problem $\left( \mathcal{P}\right) $
are (at least formally) the critical points of the Euler functional
associated to the equation of $\left( \mathcal{P}\right) $, i.e., 
\[
I\left( u\right) :=\frac{1}{2}\left\| u\right\| _{A}^{2}-\int_{\mathbb{R}%
^{N}}F\left( u\right) dx
\]
where 
\[
\left\| u\right\| _{A}^{2}:=\left( u,u\right) _{A}\quad \text{and}\quad
\left( u,v\right) _{A}:=\int_{\mathbb{R}^{N}}\left( \nabla u\cdot \nabla v+%
\frac{A}{\left| x\right| ^{\alpha }}uv\right) dx
\]
define the norm and the scalar product of the Hilbert space $H_{\alpha }^{1}$. 
More precisely, our argument is the following. We modify the function $f$
by setting $f\left( s\right) =0$ for all $s<0$, still denoting by $f$ the
modified function. This modification is not restrictive in proving Theorem 
\ref{THM: new} and will always be assumed in the rest of the paper. Then by $%
(\mathbf{h}_{0})$ there exists $M^{\prime }>0$ such that 
\[
\left| f\left( s\right) \right| \leq M\min \{\left| s\right|
^{p_{1}-1},\left| s\right| ^{p_{2}-1}\}\quad \text{and}\quad \left| F\left(
s\right) \right| \leq M^{\prime }\min \{\left| s\right| ^{p_{1}},\left|
s\right| ^{p_{2}}\}\quad \text{for all }s\in \mathbb{R},
\]
which yields in particular 
\begin{equation}
\left| f\left( s\right) \right| \leq M\left| s\right| ^{p-1}\quad \text{and}%
\quad \left| F\left( s\right) \right| \leq M^{\prime }\left| s\right|
^{p}\quad \text{for all }p\in \left[ p_{1},p_{2}\right] \text{ and }s\in %
\mathbb{R}.  \label{growth p}
\end{equation}
By the continuous embeddings $H_{\alpha }^{1}\hookrightarrow D^{1,2}(%
\mathbb{R}^{N})\hookrightarrow L^{2^{*}}(\mathbb{R}^{N})$, condition (\ref
{growth p}) with $p=2^{*}$ implies that \thinspace $I$ is of class $C^{1}$
on $H_{\alpha }^{1}$ and has Fr\'{e}chet derivative $I^{\prime }\left(
u\right) $ at any $u\in H_{\alpha }^{1}$ given by 
\[
I^{\prime }\left( u\right) v=\left( u,v\right) _{A}-\int_{\mathbb{R}%
^{N}}f\left( u\right) v\,dx\quad \text{for all }v\in H_{\alpha }^{1}.
\]
This yields that critical points of $I:H_{\alpha }^{1}\rightarrow \mathbb{R}$
satisfy (\ref{weak sol}). A standard argument shows that such critical
points are nonnegative (cf. the proof of Theorem \ref{THM: new} below) and
therefore nonzero critical points of $I$ are weak solutions to problem $%
\left( \mathcal{P}\right) $.

Then our proof proceeds essentially as follows. Given any integer $K$ such
that $1\leq K\leq N-1$, we write $x\in \mathbb{R}^{N}$ as $x=\left(
y,z\right) \in \mathbb{R}^{K}\times \mathbb{R}^{N-K}$ and in the space $%
H_{\alpha }^{1}$ we define the following closed subspaces of symmetric
functions: 
\[
H_{\mathrm{r}}:=\left\{ u\in H:u\left( x\right) =u\left( \left| x\right|
\right) \right\} \quad \text{and}\quad H_{K}:=\left\{ u\in H:u\left(
x\right) =u\left( y,z\right) =u\left( \left| y\right| ,\left| z\right|
\right) \right\} .
\]
Considering the restrictions $I_{\mid H_{\mathrm{r}}}$ and $I_{\mid H_{K}}$
with $2\leq K\leq N-2$, we prove that each of them has a nonzero
mountain-pass critical point (Section \ref{SEC: pf}), which is a weak solution 
to $\left( \mathcal{P}\right) $, since $H_{\mathrm{r}}$ and $H_{K}$ are natural constraints for $I$ thanks to the classical
Palais' Principle of Symmetric Criticality \cite{Palais}. So the proof of
Theorem \ref{THM: new} is accomplished by estimating the mountain-pass
critical levels we find for $I_{\mid H_{K}}$ and the nonzero critical
levels of $I_{\mid H_{\mathrm{r}}}$ (Sections \ref{SEC: mr} and \ref{SEC: cs}
respectively), in order to show that their sets are disjoint. This clearly
implies that the mountain-pass critical points $u_{K}$ of $I_{\mid H_{K}}$
are not radial, which also ensures that they are multiple, namely $%
u_{K_{1}}\neq u_{K_{2}}$ for $K_{1}\neq K_{2}$, since $H_{K_{1}}\cap
H_{K_{2}}=H_{\mathrm{r}}$ (see \cite[Lemma 2.1]{R}).

The above argument is essentially the same of \cite{R}. The main differences here
are due to the removal of the Ambrosetti-Rabinowitz condition $(\mathbf{h}%
_{3})$, which gives rise to the problem of the boundedness of the
Palais-Smale of the functionals $I_{\mid H_{\mathrm{r}}}$ and $I_{\mid H_{K}}
$. This will be overcome in Section \ref{SEC: pf} by applying an abstract results from \cite{GR-bddPS}
about the existence of bounded Palais-Smale sequences for $C^{1}$
functionals on Banach spaces. The weakening of assumptions $(\mathbf{h}_{1})$
and $(\mathbf{h}_{2})$ will be tackled, instead, by slightly modifying the
arguments in estimating the radial critical levels of $I$ (Section \ref{SEC: mr}) and the
cylindrical mountain-pass levels $I\left( u_{K}\right) $ (Section \ref{SEC: cs}).

{\allowdisplaybreaks
}


\section{Estimate of radial critical levels \label{SEC: mr}}

Let $N\geq 3$ and $\alpha ,A>0$. Let $f:\mathbb{R}\rightarrow \mathbb{R}$ be
continuous and satisfying $(\mathbf{h}_{0})$ and $(\mathbf{h}_{1}^{\prime })$%
.

This section is devoted to deriving an estimate from below for the nonzero
critical levels of $I_{\mid H_{\mathrm{r}}}$, namely for the value 
\[
m_{A}:=\inf \left\{ I\left( u\right) :u\in H_{\mathrm{r}}\setminus \left\{
0\right\} ,\,I^{\prime }(u)=0\right\} . 
\]

\begin{lemma}
\label{LEM: mr}Let $u\in H_{\mathrm{r}}\setminus \left\{ 0\right\} $ is a
critical point for $I$. Then $I\left( u\right) =\max_{t\geq 0}I\left(
tu\right) $.
\end{lemma}

\proof
As already observed, $u$ is nonnegative. For $t\geq 0$ define 
\[
g\left( t\right) :=I\left( tu\right) =\frac{1}{2}t^{2}\left\| u\right\|
_{A}^{2}-\int_{\mathbb{R}^{N}}F\left( tu\right) dx. 
\]
Then 
\[
g^{\prime }\left( t\right) =I^{\prime }\left( tu\right) u=t\left\| u\right\|
_{A}^{2}-\int_{\mathbb{R}^{N}}f\left( tu\right) u\,dx=t\left\| u\right\|
_{A}^{2}-\int_{\left\{ u>0\right\} }f\left( tu\right) u\,dx 
\]
for all $t\geq 0$ and $g^{\prime }\left( 1\right) =I^{\prime }\left(
u\right) u=0$. If $t\leq 1$, then $0\leq tu\leq u$ and assumption $(\mathbf{h%
}_{1}^{\prime })$ gives 
\[
g^{\prime }\left( t\right) =t\left( \left\| u\right\| _{A}^{2}-\int_{\left\{
u>0\right\} }\frac{f\left( tu\right) }{tu}u^{2}dx\right) \geq t\left(
\left\| u\right\| _{A}^{2}-\int_{\left\{ u>0\right\} }\frac{f\left( u\right) 
}{u}u^{2}dx\right) =I^{\prime }\left( u\right) u=0. 
\]
Similarly we get $g^{\prime }\left( t\right) \geq 0$ for $t\geq 1$, so that
we conclude $g\left( 1\right) =\max_{t\geq 0}g\left( t\right) $, which is
the claim.\endproof

\begin{proposition}
\label{PROP: mr>}Assume $0<\alpha <2N-2$, $\alpha \neq 2$, and let $p=\max
\{2_{\alpha }^{*},p_{1}\}$ or $p=\min \{2_{\alpha }^{*},p_{2}\}$ if $%
0<\alpha <2$ or $2<\alpha <2N-2$, respectively. Then there exists a constant 
$C_{0}>0$, independent from $A$, such that 
\[
m_{A}\geq C_{0}A^{\frac{N-2}{\alpha -2}\frac{p-2^{*}}{p-2}}.
\]
\end{proposition}

\proof
The argument is similar to the one of \cite[Proposition 3.2]{R}, except for
the conclusion, so we omit most computations. Let $u\in H_{\mathrm{r}%
}\setminus \left\{ 0\right\} $. Here $C$ will denote any positive constant
independent from $A$ and $u$. By the radial lemma \cite[Lemma 3.1]{R}, we
have 
\[
\left| u\left( x\right) \right| \leq \frac{C}{A^{1/4}}\frac{\left\|
u\right\| _{A}}{\left| x\right| ^{\frac{2N-2-\alpha }{4}}}\quad \text{almost
everywhere in }\mathbb{R}^{N} 
\]
and therefore 
\[
\int_{\mathbb{R}^{N}}\left| u\right| ^{2_{\alpha }^{*}}dx=\int_{\mathbb{R}%
^{N}}\left| u\right| ^{2_{\alpha }^{*}-2}u^{2}dx\leq C\frac{\left\|
u\right\| _{A}^{2_{\alpha }^{*}-2}}{A^{\left( 2_{\alpha }^{*}-2\right) /4}}%
\int_{\mathbb{R}^{N}}\frac{u^{2}}{\left| x\right| ^{\frac{2N-2-\alpha }{4}%
\left( 2_{\alpha }^{*}-2\right) }}dx\leq \frac{C}{A^{\frac{2N-2}{2N-2-\alpha 
}}}\left\| u\right\| _{A}^{2_{\alpha }^{*}}. 
\]
Then, both for $p=\max \{2_{\alpha }^{*},p_{1}\}<2^{*}$ and $p=\min
\{2_{\alpha }^{*},p_{2}\}>2^{*}$, we can use Sobolev inequality and argue by
interpolation. We get that there exists $\lambda \in \left[ 0,1\right) $
such that 
\[
\int_{\mathbb{R}^{N}}\left| u\right| ^{p}dx\leq \left( \int_{\mathbb{R}%
^{N}}\left| u\right| ^{2^{*}}dx\right) ^{\lambda }\left( \int_{\mathbb{R}%
^{N}}\left| u\right| ^{2_{\alpha }^{*}}dx\right) ^{1-\lambda }\leq C\frac{%
\left\| u\right\| _{A}^{p}}{A^{\frac{\left( 2N-2\right) \left( 1-\lambda
\right) }{2N-2-\alpha }}}. 
\]
Recalling condition (\ref{growth p}), this implies 
\[
\left| \int_{\mathbb{R}^{N}}F\left( u\right) dx\right| \leq M_{2}\int_{%
\mathbb{R}^{N}}\left| u\right| ^{p}dx\leq C\frac{\left\| u\right\| _{A}^{p}}{%
A^{\frac{\left( 2N-2\right) \left( 1-\lambda \right) }{2N-2-\alpha }}} 
\]
and therefore $I\left( u\right) \geq \left\| u\right\| _{A}^{2}/2-CA^{-\frac{%
\left( 2N-2\right) \left( 1-\lambda \right) }{2N-2-\alpha }}\left\|
u\right\| _{A}^{p}$. Then for every $t\geq 0$ we have 
\[
I\left( tu\right) \geq \frac{1}{2}t^{2}\left\| u\right\| _{A}^{2}-CA^{-\frac{%
\left( 2N-2\right) \left( 1-\lambda \right) }{2N-2-\alpha }}t^{p}\left\|
u\right\| _{A}^{p} 
\]
where the function of the variable $t$ on the r.h.s. has a maximum which can
be easily computed and it is given by $CA^{\frac{N-2}{\alpha -2}\frac{p-2^{*}%
}{p-2}}$. Hence, if $u$ is a critical point of $I$, Lemma \ref{LEM: mr}
gives 
\[
I\left( u\right) =\max_{t\geq 0}I\left( tu\right) \geq CA^{\frac{N-2}{\alpha
-2}\frac{p-2^{*}}{p-2}}. 
\]
This yields the conclusion, since $u\in H_{\mathrm{r}}\setminus \left\{
0\right\} $ is an arbitrary critical point and $C$ does not depend on $u$.%
\endproof

\begin{remark}
\label{RMK: p}If $p$ is as in Proposition \ref{PROP: mr>}, it is easy to
check that 
\[
\frac{N-2}{\alpha -2}\frac{p-2^{*}}{p-2}=\left\{ 
\begin{array}{lll}
\min \left\{ \frac{N-1}{\alpha },\frac{N-2}{2-\alpha }\frac{2^{*}-p_{1}}{%
p_{1}-2}\right\} \medskip  & \text{if } & 0<\alpha <2 \\ 
\min \left\{ \frac{N-1}{\alpha },\frac{N-2}{\alpha -2}\frac{p_{2}-2^{*}}{%
p_{2}-2}\right\}  & \text{if } & 2<\alpha <2N-2.
\end{array}
\right. 
\]
\end{remark}


\section{Estimate of cylindrical mountain-pass levels \label{SEC: cs}}

Let $N\geq 3$, $2\leq K\leq N-2$ and $\alpha >0$, $\alpha \neq 2$. Let $f:%
\mathbb{R}\rightarrow \mathbb{R}$ be a continuous function satisfying $(%
\mathbf{h}_{0})$ and $(\mathbf{h}_{2}^{\prime })$.

In this section we show that the functional $I_{\mid H_{K}}$ has a
mountain-pass geometry and provide an estimate from above for the
corresponding mountain-pass level.

The mountain-pass geometry of $I_{\mid H_{K}}$ near the origin is
straightforward, since condition (\ref{growth p}) with $p=2^{*}$ and the
continuous embedding $H_{K}\hookrightarrow L^{2^{*}}(\mathbb{R}^{N})$ imply
that $\exists C>0$ such that $I\left( u\right) \geq \left\| u\right\|
_{A}^{2}/2-C\left\| u\right\| _{A}^{2^{*}}$ for all $u\in H_{K}$. Hence $%
\exists R>0$ such that 
\begin{equation}
\inf_{u\in H_{K},\left\| u\right\| _{A}\leq R}I\left( u\right) =0\quad \text{%
and\quad }\inf_{u\in H_{K},\left\| u\right\| _{A}=R}I\left( u\right) >0.
\label{mp0}
\end{equation}

Now we define a suitable $\overline{u}_{K}\in H_{K}$ such that $I\left( 
\overline{u}_{K}\right) <0$. Then one has $\left\| \overline{u}_{K}\right\|
_{A}>R$ by \eqref{mp0}, so that $I_{\mid H_{K}}$ has the mountain-pass
geometry and we can define the mountain-pass level 
\begin{equation}
c_{A,K}:=\inf_{\gamma \in \Gamma }\max_{t\in \left[ 0,1\right] }I\left(
\gamma \left( t\right) \right) >0\quad \text{where\quad }\Gamma :=\left\{
\gamma \in C\left( \left[ 0,1\right] ;H_{K}\right) :\gamma \left( 0\right)
=0,\,\gamma \left( 1\right) =\overline{u}_{K}\right\} .  \label{cs:=}
\end{equation}

The definition of $\overline{u}_{K}$ is inspired by some arguments of \cite
{Bad-Serra-mult}. Denote by $\phi :D\rightarrow \mathbb{R}^{2}\setminus
\left\{ 0\right\} $ the change to polar coordinates in $\mathbb{R}%
^{2}\setminus \left\{ 0\right\} $, namely $\phi \left( \rho ,\theta \right)
=\left( \rho \cos \theta ,\rho \sin \theta \right) $ for all $\left( \rho
,\theta \right) \in D:=\left( 0,+\infty \right) \times \left[ 0,2\pi \right) 
$. Define $E:=\left( 1/4,3/4\right) \times \left( \pi /6,\pi /3\right) $ and
take any $\psi :\mathbb{R}^{2}\rightarrow \mathbb{R}$ such that $\psi \in
C_{c}^{\infty }\left( E\right) $, $\psi \neq 0$ and $0\leq \psi <s_{*}$,
where $s_{*}$ is given in assumption $(\mathbf{h}_{2}^{\prime })$. For every 
$A>1$, define 
\[
E_{A}:=\left\{ \left( \rho ,\theta \right) \in \mathbb{R}^{2}:\left( \frac{1%
}{4}\right) ^{1/\sqrt{A}}<\rho <\left( \frac{3}{4}\right) ^{1/\sqrt{A}},%
\frac{\pi }{6\sqrt{A}}<\theta <\frac{\pi }{3\sqrt{A}}\right\} 
\]
and a function $\psi _{A}\in C_{c}^{\infty }\left( E_{A}\right) $ by setting 
\[
\psi _{A}\left( \rho ,\theta \right) :=\psi \left( \rho ^{\sqrt{A}},\theta 
\sqrt{A}\right) . 
\]
Finally define 
\[
v_{A}\left( y,z\right) :=\psi _{A}\left( \phi ^{-1}\left( \left| y\right|
,\left| z\right| \right) \right) \quad \text{for }x=\left( y,z\right) \in (%
\mathbb{R}^{K}\times \mathbb{R}^{N-K})\setminus \left\{ 0\right\} ,\quad
v_{A}\left( 0\right) :=0. 
\]
Then $v_{A}\in C_{c}^{\infty }\left( \Omega _{A}\right) \cap H_{K}$, where $%
\Omega _{A}:=\left\{ \left( y,z\right) \in \mathbb{R}^{K}\times \mathbb{R}%
^{N-K}:\left( \left| y\right| ,\left| z\right| \right) \in \phi \left(
E_{A}\right) \right\} $. Notice that $v_{A}\neq 0$ and $0\leq v_{A}<s_{*}$.

By using spherical coordinates in $\mathbb{R}^{K}$ and $\mathbb{R}^{N-K}$,
and then making the change of variables $\left( r,\varphi \right) =(\rho ^{%
\sqrt{A}},\theta \sqrt{A})$, one computes 
\begin{eqnarray}
\int_{\mathbb{R}^{N}}\frac{v_{A}^{2}}{\left| x\right| ^{\alpha }}dx
&=&\sigma _{K}\sigma _{N-K}\int_{E_{A}}\frac{\psi \left( \rho ^{\sqrt{A}%
},\theta \sqrt{A}\right) ^{2}}{\rho ^{\alpha -N+1}}H\left( \theta \right)
d\rho \,d\theta  \nonumber \\
&=&\frac{\sigma _{K}\sigma _{N-K}}{A}\int_{E}\frac{\psi \left( r,\varphi
\right) ^{2}}{r^{\left( \alpha -N\right) /A^{1/2}+1}}H(A^{-1/2}\varphi
)dr\,d\varphi  \label{comp1}
\end{eqnarray}
where $H\left( \theta \right) :=\left( \cos \theta \,\right) ^{K-1}\left(
\sin \theta \right) ^{N-K-1}$ (see \cite{R} for more detailed computation).
Here and in the following $\sigma _{d}$ denotes the $\left( d-1\right) $%
-dimensional measure of the unit sphere of $\mathbb{R}^{d}$. Similarly one
obtains 
\begin{equation}
\int_{\mathbb{R}^{N}}F\left( v_{A}\right) dx=\frac{\sigma _{K}\sigma _{N-K}}{%
A}\int_{E}F\left( \psi \left( r,\varphi \right) \right)
r^{N/A^{1/2}-1}H\left( A^{-1/2}\varphi \right) dr\,d\varphi  \label{comp2}
\end{equation}
and 
\begin{equation}
\int_{\mathbb{R}^{N}}\left| \nabla v_{A}\right| ^{2}dx=\sigma _{K}\sigma
_{N-K}\int_{E}\left( \psi _{r}\left( r,\varphi \right) ^{2}+\frac{1}{r^{2}}%
\psi _{\varphi }\left( r,\varphi \right) ^{2}\right)
r^{(N-2)/A^{1/2}+1}H\left( A^{-1/2}\varphi \right) dr\,d\varphi
\label{comp3}
\end{equation}
where we denote $\psi _{r}=\frac{\partial \psi }{\partial r}$ and $\psi
_{\varphi }=\frac{\partial \psi }{\partial \varphi }$ for brevity.

Then we have 
\[
\frac{\left\| v_{A}\right\| _{A}^{2}}{\int_{\mathbb{R}^{N}}F\left(
v_{A}\right) dx}=A\frac{\int_{E}\left( \left( \psi _{r}^{2}+\frac{1}{r^{2}}%
\psi _{\varphi }^{2}\right) r^{(N-2)/A^{1/2}+1}+\psi ^{2}r^{\left( N-\alpha
\right) /A^{1/2}-1}\right) H(A^{-1/2}\varphi )dr\,d\varphi }{\int_{E}F\left(
\psi \right) r^{N/A^{1/2}-1}H\left( A^{-1/2}\varphi \right) dr\,d\varphi }. 
\]
As in the integration set $E$ one has $A^{-1/2}\pi /6<A^{-1/2}\varphi
<A^{-1/2}\pi /3$, for $A>1$ large enough we have that $A^{-1/2}\varphi
/2<\sin (A^{-1/2}\varphi )<A^{-1/2}\varphi $ and $1/2<(\cos A^{-1/2}\varphi
)<1$, and thus there exist two constants $\overline{C}_{1},\overline{C}%
_{2}>0 $ such that 
\[
\overline{C}_{1}A^{-\frac{N-K-1}{2}}<H\left( A^{-1/2}\varphi \right) <%
\overline{C}_{2}A^{-\frac{N-K-1}{2}}. 
\]
Similarly, since $1/4<r<3/4$ in $E$, the terms $r^{(N-2)A^{-1/2}+1}$, $%
r^{\left( N-\alpha \right) A^{-1/2}-1}$ and $r^{NA^{-1/2}-1}$ are bounded
and bounded away from zero by positive constants independent of $A>1$, say $%
\overline{C}_{3}$ and $\overline{C}_{4}$ respectively, so that we conclude 
\[
\frac{\left\| v_{A}\right\| _{A}^{2}}{\int_{\mathbb{R}^{N}}F\left(
v_{A}\right) dx}\geq A\frac{\overline{C}_{1}\overline{C}_{4}\int_{E}\left(
\psi _{r}^{2}+\frac{1}{r^{2}}\psi _{\varphi }^{2}+\psi ^{2}\right)
dr\,d\varphi }{\overline{C}_{2}\overline{C}_{3}\int_{E}F\left( \psi \right)
dr\,d\varphi }\rightarrow +\infty \quad \text{as }A\rightarrow +\infty . 
\]
This allows us to fix, hereafter, a threshold $A_{K}>0$ such that 
\begin{equation}
\frac{\left\| v_{A}\right\| _{A}^{2}}{\int_{\mathbb{R}^{N}}F\left(
v_{A}\right) dx}>1\quad \text{for every }A>A_{K}.  \label{def A0}
\end{equation}

\begin{proposition}
\label{PROP: cs <2}Assume $0<\alpha <2$ and let $A>A_{K}$. Define $\overline{%
u}_{K}\in H_{K}$ by setting 
\[
\overline{u}_{K}\left( x\right) :=v_{A}\left( \frac{x}{\lambda }\right)
\quad \text{with\quad }\lambda =\frac{\left\| v_{A}\right\| _{A}^{2/\alpha }%
}{\left( \int_{\mathbb{R}^{N}}F\left( v_{A}\right) dx\right) ^{1/\alpha }}.
\]
Then $I\left( \overline{u}_{K}\right) <0$ and the corresponding
mountain-pass level (\ref{cs:=}) satisfies 
\[
c_{A,K}\leq C_{1}A^{\frac{K-1}{2}+N\left( \frac{1}{\alpha }-\frac{1}{2}%
\right) }
\]
where the constant $C_{1}>0$ does not depend on $A$.
\end{proposition}

\proof
As $\lambda >1$, an obvious change of variables yields 
\begin{eqnarray*}
I\left( \overline{u}_{K}\right) &=&\frac{\lambda ^{N-2}}{2}\int_{\mathbb{R}%
^{N}}\left| \nabla v_{A}\right| ^{2}dx+\frac{\lambda ^{N-\alpha }}{2}\int_{%
\mathbb{R}^{N}}\frac{A}{\left| x\right| ^{\alpha }}v_{A}^{2}dx-\lambda
^{N}\int_{\mathbb{R}^{N}}F\left( v_{A}\right) dx \\
&\leq &\frac{\lambda ^{N-\alpha }}{2}\left( \int_{\mathbb{R}^{N}}\left|
\nabla v_{A}\right| ^{2}dx+\int_{\mathbb{R}^{N}}\frac{A}{\left| x\right|
^{\alpha }}v_{A}^{2}dx\right) -\lambda ^{N}\int_{\mathbb{R}^{N}}F\left(
v_{A}\right) dx \\
&=&\frac{\lambda ^{N}}{2}\left( \lambda ^{-\alpha }\left\| v_{A}\right\|
_{A}^{2}-2\int_{\mathbb{R}^{N}}F\left( v_{A}\right) dx\right) =-\frac{%
\lambda ^{N}}{2}\int_{\mathbb{R}^{N}}F\left( v_{A}\right) dx<0.
\end{eqnarray*}
In order to estimate $c_{A,K}$, observe that hypothesis $(\mathbf{h}%
_{2}^{\prime })$ implies $F\left( ts\right) \geq t^{\mu }F\left( s\right) $
for all $0\leq s<s_{*}$ and $t\in \left[ 0,1\right] $. Then consider the
straight path $\gamma \left( t\right) :=t\overline{u}_{K}$, $t\in \left[
0,1\right] $. Since $0\leq \overline{u}_{K}<s_{*}$ (recall that $0\leq
v_{A}<s_{*}$), we have 
\begin{equation}
I\left( \gamma \left( t\right) \right) =\frac{1}{2}t^{2}\left\| \overline{u}%
_{K}\right\| _{A}^{2}-\int_{\mathbb{R}^{N}}F\left( t\overline{u}_{K}\right)
dx\leq \frac{1}{2}t^{2}\left\| \overline{u}_{K}\right\| _{A}^{2}-t^{\mu
}\int_{\mathbb{R}^{N}}F\left( \overline{u}_{K}\right) dx.  \label{I(gamma)}
\end{equation}
Note that $F\left( \overline{u}_{K}\right) \neq 0$ since $\overline{u}%
_{K}\neq 0$. The function of the variable $t$ on the r.h.s of (\ref{I(gamma)}%
) has maximum 
\[
m\left\| \overline{u}_{K}\right\| _{A}^{2}\left( \frac{\left\| \overline{u}%
_{K}\right\| _{A}^{2}}{\int_{\mathbb{R}^{N}}F\left( \overline{u}_{K}\right)
dx}\right) ^{2/(\mu -2)} 
\]
where $m:=\left( 1/\mu \right) ^{2/(\mu -2)}\left( 1/2-1/\mu \right) $ for
brevity. Therefore, recalling that $\lambda >1$, we obtain 
\begin{eqnarray*}
c_{A,K} &\leq &\max_{t\in \left[ 0,1\right] }I\left( \gamma \left( t\right)
\right) \leq m\frac{\left\| \overline{u}_{K}\right\| _{A}^{2\mu /(\mu -2)}}{%
\left( \int_{\mathbb{R}^{N}}F\left( \overline{u}_{K}\right) dx\right)
^{2/(\mu -2)}} \\
&\leq &m\frac{\lambda ^{\mu \left( N-\alpha \right) /(\mu -2)}\left( \int_{%
\mathbb{R}^{N}}\left| \nabla v_{A}\right| ^{2}dx+\int_{\mathbb{R}%
^{N}}A\left| x\right| ^{-\alpha }v_{A}^{2}dx\right) ^{\mu /(\mu -2)}}{%
\lambda ^{2N/(\mu -2)}\left( \int_{\mathbb{R}^{N}}F\left( v_{A}\right)
dx\right) ^{2/(\mu -2)}} \\
&=&m\lambda ^{\frac{\mu \left( N-\alpha \right) -2N}{\mu -2}}\frac{\left\|
v_{A}\right\| _{A}^{2\mu /(\mu -2)}}{\left( \int_{\mathbb{R}^{N}}F\left(
v_{A}\right) dx\right) ^{2/(\mu -2)}}.
\end{eqnarray*}
Inserting the definition of $\lambda $ and using computations (\ref{comp1})-(%
\ref{comp3}), we get 
\begin{eqnarray*}
c_{A,K} &\leq &m\frac{\left\| v_{A}\right\| _{A}^{\frac{2N}{\alpha }}}{%
\left( \int_{\mathbb{R}^{N}}F\left( v_{A}\right) dx\right) ^{\frac{N-\alpha 
}{\alpha }}} \\
&=&m\sigma _{K}\sigma _{N-K}\frac{\left( \int_{E}\left( \left( \psi _{r}^{2}+%
\frac{1}{r^{2}}\psi _{\varphi }^{2}\right) r^{(N-2)A^{-1/2}+1}+\psi
^{2}r^{\left( N-\alpha \right) A^{-1/2}-1}\right) H\left( A^{-1/2}\varphi
\right) dr\,d\varphi \right) ^{\frac{N}{\alpha }}}{A^{^{-\frac{N-\alpha }{%
\alpha }}}\left( \int_{E}F\left( \psi \right) r^{NA^{-1/2}-1}H\left(
A^{-1/2}\varphi \right) dr\,d\varphi \right) ^{\frac{N-\alpha }{\alpha }}}.
\end{eqnarray*}
Finally, taking as above the four constants $\overline{C}_{1},...,\overline{C%
}_{4}>0$ independent of $A$, we conclude 
\begin{eqnarray*}
c_{A,K} &\leq &m\sigma _{K}\sigma _{N-K}\frac{\left( \overline{C}_{2}%
\overline{C}_{3}\int_{E}\left( \left( \psi _{r}^{2}+\frac{1}{r^{2}}\psi
_{\varphi }^{2}\right) +\psi ^{2}r\right) A^{-\frac{N-K-1}{2}}dr\,d\varphi
\right) ^{\frac{N}{\alpha }}}{A^{-\frac{N-\alpha }{\alpha }}\left( \overline{%
C}_{1}\overline{C}_{4}\int_{E}F\left( \psi \right) A^{-\frac{N-K-1}{2}%
}dr\,d\varphi \right) ^{\frac{N-\alpha }{\alpha }}} \\
&=&CA^{\frac{K-1}{2}+N\left( \frac{1}{\alpha }-\frac{1}{2}\right) }\frac{%
\left( \int_{E}\left( \left( \psi _{r}^{2}+\frac{1}{r^{2}}\psi _{\varphi
}^{2}\right) +\psi ^{2}r\right) dr\,d\varphi \right) ^{\frac{N}{\alpha }}}{%
\left( \int_{E}F\left( \psi \right) dr\,d\varphi \right) ^{\frac{N-\alpha }{%
\alpha }}}
\end{eqnarray*}
with obvious definition of the constant $C$. As the last ratio does not
depend on $A$, the conclusion ensues. \endproof

\begin{proposition}
\label{PROP: cs >2}Assume $\alpha >2$ and let $A>A_{K}$. Define $\overline{u}%
\in H_{K}$ by setting 
\[
\overline{u}_{K}\left( x\right) :=w_{A}\left( \frac{x}{\lambda }\right)
\quad \text{with\quad }\lambda =\frac{\left\| v_{A}\right\| _{A}}{\left(
\int_{\mathbb{R}^{N}}F\left( v_{A}\right) dx\right) ^{1/2}}.
\]
Then $I\left( \overline{u}_{K}\right) <0$ and the corresponding
mountain-pass level (\ref{cs:=}) satisfies 
\[
c_{A,K}\leq C_{2}A^{\frac{K-1}{2}}
\]
where the constant $C_{2}>0$ does not depend on $A$.
\end{proposition}

\proof
The proof is very similar to the one of Proposition \ref{PROP: cs <2} above
(see also \cite[Proposition 4.3]{R}), so we leave it to the reader for
brevity.\endproof


\section{Proof of Theorem \ref{THM: main} \label{SEC: pf}}

This section is devoted to the proof of Theorem \ref{THM: main}, so assume
all the hypotheses of the theorem.

Let $K$ be any integer such that $2\leq K\leq N-2$. Assume $A>A_{K}$ (where $%
A_{K}$ is defined by (\ref{def A0})) and consider the mountain-pass level $%
c_{A,K}$ defined by (\ref{cs:=}), with $\overline{u}_{K}\in H_{K}$ given by
Proposition \ref{PROP: cs <2} or \ref{PROP: cs >2} if $\alpha \in \left(
2/(N-1),2\right) $ or $\alpha \in \left( 2,2N-2\right) $, respectively.

We first show that $c_{A,K}$ is a critical level for the energy functional $%
I $. To do this, we will make use of the sum space 
\[
L^{p_{1}}+L^{p_{2}}:=\left\{ u_{1}+u_{2}:u_{1}\in L^{p_{1}}\left( \mathbb{R}%
^{N}\right) ,\,u_{2}\in L^{p_{2}}\left( \mathbb{R}^{N}\right) \right\} . 
\]
We recall that it is a Banach space with respect to the norm 
\[
\left\| u\right\| _{L^{p_{1}}+L^{p_{2}}}:=\inf_{u_{1}+u_{2}=u}\max \left\{
\left\| u_{1}\right\| _{L^{p_{1}}(\mathbb{R}^{N})},\left\| u_{2}\right\|
_{L^{p_{2}}(\mathbb{R}^{N})}\right\} 
\]
(see \cite[Corollary 2.11]{BPR}) and that the continuous embedding $L^{p}(%
\mathbb{R}^{N})\hookrightarrow L^{p_{1}}+L^{p_{2}}$ holds for all $p\in
\left[ p_{1},p_{2}\right] $ (see \cite[Proposition 2.17]{BPR}), in
particular for $p=2^{*}$. Moreover, for every $u\in L^{p_{1}}+L^{p_{2}}$ and
every $\varphi \in L^{p_{1}^{\prime }}(\mathbb{R}^{N})\cap L^{p_{2}^{\prime
}}(\mathbb{R}^{N})$ one has 
\begin{equation}
\int_{\mathbb{R}^{N}}\left| u\varphi \right| dx\leq \left\| u\right\|
_{L^{p_{1}}+L^{p_{2}}}\left( \left\| \varphi \right\| _{L^{p_{1}^{\prime }}(%
\mathbb{R}^{N})}+\left\| \varphi \right\| _{L^{p_{2}^{\prime }}(\mathbb{R}%
^{N})}\right)  \label{dual}
\end{equation}
where $p_{i}^{\prime }=p_{i}/(p_{i}-1)$ is the H\"{o}lder conjugate exponent
of $p_{i}$ (see \cite[Lemma 2.9]{BPR}).

We will also need the following abstract result from \cite{GR-bddPS},
concerning the existence of bounded Palais-Smale sequences for $C^{1}$
functionals on Banach spaces.

\begin{lemma}
\label{LEM: bdd PS}Let $J:X\rightarrow \mathbb{R}$ be a $C^{1}$ functional
on a Banach space $\left( X,\left\| \cdot \right\| \right) $, having the
form 
\[
J\left( u\right) =\frac{1}{q}\left\| u\right\| ^{q}-B\left( u\right) 
\]
for some $q>0$ and some continuous functional $B:X\rightarrow \mathbb{R}$.
Assume that there exists a sequence of continuous mappings $\psi
_{n}:X\rightarrow X$ such that $\forall n$ there exist $\alpha _{n}>\beta
_{n}>0$ satisfying 
\begin{equation}
\left\| u\right\| ^{q}\geq \alpha _{n}\left\| \psi _{n}\left( u\right)
\right\| ^{q}\text{\quad and\quad }B\left( u\right) \leq \beta _{n}B\left(
\psi _{n}\left( u\right) \right) \text{\quad for all }u\in X
\label{alfabeta}
\end{equation}
and 
\[
\lim\limits_{n\rightarrow \infty }\alpha _{n}=\lim\limits_{n\rightarrow
\infty }\beta _{n}=1\,,\qquad \liminf\limits_{n\rightarrow \infty }\frac{%
\left| 1-\beta _{n}\right| }{\alpha _{n}-\beta _{n}}<\infty .
\]
If there exist $r>0$ and $\overline{u}\in X$ with $\left\| \overline{u}%
\right\| >r$ such that 
\[
\inf_{\left\| u\right\| =r}J\left( u\right) >J\left( 0\right) \geq J\left( 
\overline{u}\right) \quad \text{and}\quad \lim\limits_{n\rightarrow \infty
}\left\| \psi _{n}\left( 0\right) \right\| =\lim\limits_{n\rightarrow \infty
}\left\| \psi _{n}\left( \overline{u}\right) -\overline{u}\right\| =0,
\]
then $J$ has a bounded Palais-Smale sequence $\left\{ w_{n}\right\} \subset X
$ at level 
\[
c_{J,\bar{u}}=\inf\limits_{\gamma \in \Gamma }\max\limits_{u\in \gamma
\left( \left[ 0,1\right] \right) }J\left( u\right) ,~~\Gamma =\left\{ \gamma
\in C\left( \left[ 0,1\right] ,X\right) :\gamma \left( 0\right) =0,\,\gamma
\left( 1\right) =\bar{u}\right\} .
\]
\end{lemma}

\proof%
It is a particular case of \cite[Theorem 1.1]{GR-bddPS}.%
\endproof%

\begin{lemma}
\label{LEM: ck critical}$c_{A,K}$ is a critical level for the functional $%
I_{\mid H_{K}}$.
\end{lemma}

\proof
We want to apply Lemma \ref{LEM: bdd PS} to the functional $I_{\mid H_{K}}$
with $\overline{u}=\overline{u}_{K}$ (and of course $q=2$ and $B\left(
u\right) =\int_{\mathbb{R}^{N}}F\left( u\right) dx$). Take any real sequence 
$\left\{ t_{n}\right\} $ such that $t_{n}\rightarrow 1$ and $t_{n}>1$, and
for every $u\in H_{K}$ define $\psi _{n}\left( u\right) \in H_{K}$ by
setting $\psi _{n}\left( u\right) \left( x\right) :=u(x/t_{n})$. We have 
\begin{equation}
\left\| \psi _{n}\left( u\right) \right\| _{A}^{2}=t_{n}^{N-2}\int_{%
\mathbb{R}^{N}}\left| \nabla u\right| ^{2}dx+t_{n}^{N-\alpha }\int_{%
\mathbb{R}^{N}}\frac{A}{\left| x\right| ^{\alpha }}u^{2}dx\leq t_{n}^{\max
\left\{ N-2,N-\alpha \right\} }\left\| u\right\| _{A}^{2}  \label{A}
\end{equation}
and 
\[
\int_{\mathbb{R}^{N}}F\left( \psi _{n}\left( u\right) \right)
dx=t_{n}^{N}\int_{\mathbb{R}^{N}}F\left( u\right) dx, 
\]
so that (\ref{alfabeta}) holds with $\alpha _{n}=t_{n}^{-\max \left\{
N-2,N-\alpha \right\} }$ and $\beta _{n}=t_{n}^{-N}$. Note that (\ref{A})
also ensures that the linear mapping $\psi _{n}:H_{K}\rightarrow H_{K}$ is
continuous. Since $t_{n}>1$ and $\max \left\{ N-2,N-\alpha \right\} <N$, we
have $\alpha _{n}>\beta _{n}$. Moreover $\lim\limits_{n\rightarrow \infty
}\alpha _{n}=\lim\limits_{n\rightarrow \infty }\beta _{n}=1$ and 
\[
\lim\limits_{n\rightarrow \infty }\frac{\left| 1-\beta _{n}\right| }{\alpha
_{n}-\beta _{n}}=\lim\limits_{n\rightarrow \infty }\frac{1-t_{n}^{-N}}{%
t_{n}^{-\max \left\{ N-2,N-\alpha \right\} }-t_{n}^{-N}}=\frac{N}{N-\max
\left\{ N-2,N-\alpha \right\} }<\infty . 
\]
Recalling (\ref{mp0}) and the fact that $I\left( \overline{u}_{K}\right) <0$%
, it remains to check that $\lim\limits_{n\rightarrow \infty }\left\| \psi
_{n}\left( \overline{u}_{K}\right) -\overline{u}_{K}\right\| _{A}=0$. By
definition of $\overline{u}_{K}$, there exist $l_{1},l_{2}>0$ such that both 
$\overline{u}_{K}\in C_{c}^{\infty }(\Omega _{l_{1},l_{2}})$ where $\Omega
_{l_{1},l_{2}}:=\left\{ x\in \mathbb{R}^{N}:l_{1}<\left| x\right|
<l_{2}\right\} $. Then both $\overline{u}_{K}$ and $\psi _{n}\left( 
\overline{u}_{K}\right) $ belong to $C_{c}^{\infty }(\Omega
_{l_{1}/2,2l_{2}})$ for $n$ sufficiently large, and thus we get 
\[
\int_{\mathbb{R}^{N}}\frac{\left( \psi _{n}\left( \overline{u}_{K}\right) -%
\overline{u}_{K}\right) ^{2}}{\left| x\right| ^{\alpha }}dx=\int_{\Omega
_{l_{1}/2,2l_{2}}}\frac{\left( \psi _{n}\left( \overline{u}_{K}\right) -%
\overline{u}_{K}\right) ^{2}}{\left| x\right| ^{\alpha }}dx\leq \frac{%
2^{\alpha }}{l_{1}^{\alpha }}\int_{\Omega _{l_{1}/2,2l_{2}}}\left( \psi
_{n}\left( \overline{u}_{K}\right) -\overline{u}_{K}\right)
^{2}dx\rightarrow 0 
\]
and 
\[
\int_{\mathbb{R}^{N}}\left| \nabla \psi _{n}\left( \overline{u}_{K}\right)
-\nabla \overline{u}_{K}\right| ^{2}dx=\int_{\Omega _{l_{1}/2,2l_{2}}}\left|
\nabla \psi _{n}\left( \overline{u}_{K}\right) -\nabla \overline{u}%
_{K}\right| ^{2}dx\rightarrow 0 
\]
as $n\rightarrow \infty $, since $\psi _{n}\left( \overline{u}_{K}\right) =%
\overline{u}_{K}\left( t_{n}^{-1}\cdot \right) \rightarrow \overline{u}_{K}$
e $\nabla \left( \psi _{n}\left( \overline{u}_{K}\right) \right)
=t_{n}^{-1}\nabla \overline{u}_{K}\left( t_{n}^{-1}\cdot \right) \rightarrow
\nabla \overline{u}_{K}$ in $L^{\infty }(\mathbb{R}^{N})$. So we can apply
apply Lemma \ref{LEM: bdd PS} and deduce that there exists a sequence $%
\left\{ u_{n}\right\} \subset H_{K}$ such that $\left\{ \left\|
u_{n}\right\| _{A}\right\} $ is bounded, $I\left( u_{n}\right) \rightarrow
c_{A,K}$ and $I^{\prime }\left( u_{n}\right) \rightarrow 0$ in the dual
space of $H_{K}$. Now we exploit the fact that the space $H_{K}$ is
compactly embedded into $L^{p_{1}}+L^{p_{2}}$, since $p_{1}<2^{*}<p_{2}$ and
so is the subspace of $D^{1,2}(\mathbb{R}^{N})$ made up of the mappings with
the same symmetries of $H_{K}$ (see \cite[Theorem A.1]{Azz-Pomp}). Hence
there exists $u\in H_{K}$ such that (up to a subsequence) $%
u_{n}\rightharpoonup u$ in $H_{K}$ and $u_{n}\rightarrow u$ in $%
L^{p_{1}}+L^{p_{2}}$. This implies that $\left\{ f\left( u_{n}\right)
\right\} $ is bounded in both $L^{p_{1}^{\prime }}(\mathbb{R}^{N})$ and $%
L^{p_{2}^{\prime }}(\mathbb{R}^{N})$, since assumption $(\mathbf{h}_{0})$
ensures that the Nemytski\u{\i} operator $v\mapsto f\left( v\right) $ is
continuous from $L^{p_{1}}+L^{p_{2}}$ into $L^{p_{1}^{\prime }}(\mathbb{R}%
^{N})\cap L^{p_{2}^{\prime }}(\mathbb{R}^{N})$ (see \cite[Corollary 3.7]{BPR}%
). Then by (\ref{dual}) we get 
\begin{eqnarray*}
\left| \int_{\mathbb{R}^{N}}f\left( u_{n}\right) \left( u_{n}-u\right)
dx\right| &\leq &\left\| u_{n}-u\right\| _{L^{p_{1}}+L^{p_{2}}}\left(
\left\| f\left( u_{n}\right) \right\| _{L^{p_{1}^{\prime }}(\mathbb{R}%
^{N})}+\left\| f\left( u_{n}\right) \right\| _{L^{p_{2}^{\prime }}(\mathbb{R}%
^{N})}\right) \\
&\leq &\left( \mathrm{const.}\right) \left\| u_{n}-u\right\|
_{L^{p_{1}}+L^{p_{2}}}\rightarrow 0
\end{eqnarray*}
and therefore 
\begin{eqnarray*}
\left\| u_{n}-u\right\| _{A}^{2} &=&\left( u_{n},u_{n}-u\right) _{A}-\left(
u,u_{n}-u\right) _{A} \\
&=&I^{\prime }\left( u_{n}\right) \left( u_{n}-u\right) +\int_{\mathbb{R}%
^{N}}f\left( u_{n}\right) \left( u_{n}-u\right) dx-\left( u,u_{n}-u\right)
_{A}\rightarrow 0,
\end{eqnarray*}
where $\left( u,u_{n}-u\right) _{A}\rightarrow 0$ since $u_{n}%
\rightharpoonup u$ in $H_{K}$, and $I^{\prime }\left( u_{n}\right) \left(
u_{n}-u\right) \rightarrow 0$ because $I^{\prime }\left( u_{n}\right)
\rightarrow 0$ in the dual space of $H_{K}$ and $\left\{ u_{n}-u\right\} $
is bounded in $H_{K}$. This implies that $u_{n}\rightarrow u$ in $H_{K}$ and
thus concludes the proof.\endproof

Now, using again Lemma \ref{LEM: bdd PS} and the compactness results of \cite
{Su-Wang-Will-p}, we show that $I$ also has a radial critical point.

\begin{lemma}
\label{LEM: r critical}The functional $I_{\mid H_{\mathrm{r}}}$ has at least
a nonzero critical point.
\end{lemma}

\proof%
The arguments are absolutely similar to the ones used in proving Lemma \ref
{LEM: ck critical}, so we will be very sketchy. By condition (\ref{growth p}%
) with $p=2^{*}$ and the continuous embedding $H_{\mathrm{r}}\hookrightarrow
L^{2^{*}}(\mathbb{R}^{N})$, we readily have that $\exists R>0$ such that $%
\inf_{\left\| u\right\| _{A}=R}I_{\mid H_{\mathrm{r}}}\left( u\right) >0$
and $I_{\mid H_{\mathrm{r}}}\left( u\right) \geq 0$ for $\left\| u\right\|
_{A}\leq R$. For every $u\in H_{\mathrm{r}}$ and $t>1$ define $\psi
_{t}\left( u\right) \left( x\right) :=u(x/t)$. Take $u_{0}\in C_{c}^{\infty
}(\mathbb{R}^{N}\setminus \left\{ 0\right\} )\cap H_{\mathrm{r}}$ such that $%
u_{0}\geq 0$ and $u_{0}\neq 0$. We have 
\[
I\left( \psi _{t}\left( u_{0}\right) \right) \leq \frac{1}{2}t^{\max \left\{
N-2,N-\alpha \right\} }\left\| u_{0}\right\| _{A}^{2}-t^{N}\int_{\mathbb{R}%
^{N}}F\left( u\right) dx 
\]
and therefore $I\left( \overline{u}\right) <0$ for $\overline{u}:=\psi
_{t}\left( u_{0}\right) $ with $t$ large enough. Then, letting $\left\{
t_{n}\right\} $ be any real sequence such that $t_{n}\rightarrow 1$ and $%
t_{n}>1$, as in the proof of Lemma \ref{LEM: ck critical} one checks that
Lemma \ref{LEM: bdd PS} applies with $\psi _{n}=\psi _{t_{n}}$, $\alpha
_{n}=t_{n}^{-\max \left\{ N-2,N-\alpha \right\} }$ and $\beta
_{n}=t_{n}^{-N} $. So $I_{\mid H_{\mathrm{r}}}$ has a Palais-Smale sequence $%
\left\{ u_{n}\right\} $ at a nonzero level. Now observe that assumption $(%
\mathbf{h}_{0})$ with $p_{1},p_{2}$ as in the theorem ensures that one can
find $p\in [p_{1},p_{2}]$ such that $\left| f\left( u\right) \right| \leq
\left( \mathrm{const.}\right) u^{p-1}$ (recall (\ref{growth p})) and (\ref
{sww}) holds, in such a way that the embedding $H_{\mathrm{r}%
}\hookrightarrow L^{p}(\mathbb{R}^{N})$ is compact by the compactness
results of \cite{Su-Wang-Will-p}. Hence it is a standard exercise to
conclude that, up to a subsequence, $u_{n}$ converges to a nonzero critical
point of $I_{\mid H_{\mathrm{r}}}$.\endproof

\proof[Proof of Theorem \ref{THM: main}] The proof is essentially the same
of Theorem 1.1 of \cite{R}; we repeat it here just for the sake of
completeness. On the one hand, the restriction $I_{\mid H_{\mathrm{r}}}$ has
a critical point $u_{\mathrm{r}}\neq 0$ by Lemma \ref{LEM: r critical}. On
the other hand, one checks that for $\alpha \in \left( \frac{2}{N-1}%
,2N-2\right) \setminus \left\{ 2\right\} $ the integer $\nu $ defined in (%
\ref{n_alfa}) is at least $1$ (see \cite[Lemma 5.2]{R}), so that there are $%
\nu $ integers $K$ (precisely $K=2,...,\nu +1$) such that 
\[
\frac{K-1}{2}+N\left( \frac{1}{\alpha }-\frac{1}{2}\right) <\min \left\{ 
\frac{N-1}{\alpha },\frac{N-2}{2-\alpha }\frac{2^{*}-p_{1}}{p_{1}-2}\right\}
\quad \text{if\quad }\frac{2}{N-1}<\alpha <2 
\]
and 
\[
\frac{K-1}{2}<\min \left\{ \frac{N-1}{\alpha },\frac{N-2}{\alpha -2}\frac{%
p_{2}-2^{*}}{p_{2}-2}\right\} \quad \text{if\quad }2<\alpha <2N-2. 
\]
Let $K$ be any of such integers. By Lemma \ref{LEM: ck critical}, there
exists $u_{K}\in H_{K}\setminus \left\{ 0\right\} $ such that $I\left(
u_{K}\right) =c_{A,K}>0$ and $I_{\mid H_{K}}^{\prime }\left( u_{K}\right) =0$%
. Both $u_{\mathrm{r}}$ and $u_{K}$ are also critical points for the
functional $I:H_{\alpha }^{1}\rightarrow \mathbb{R}$, by the Palais'
Principle of Symmetric Criticality \cite{Palais}. Moreover, it easy to check
that they are nonnegative: test $I^{\prime }\left( u_{K}\right) $ with the
negative part $u_{K}^{-}\in H_{\alpha }^{1}$ of $u_{K}$ and use the fact
that $f\left( s\right) =0$ for $s<0$ to get $I^{\prime }\left( u_{K}\right)
u_{K}^{-}=-\left\| u_{K}^{-}\right\| _{A}^{2}=0$; the same for $u_{\mathrm{r}%
}$. Therefore $u_{\mathrm{r}}$ and $u_{K}$ are weak solutions to problem $%
\left( \mathcal{P}\right) $. By Remark \ref{RMK: p} and Propositions \ref
{PROP: mr>}, \ref{PROP: cs <2} and \ref{PROP: cs >2}, there exists $%
\widetilde{A}_{K}>A_{K}$ such that $c_{A,K}<m_{A}$ for every $A>\widetilde{A}%
_{K}$. Now, setting $A_{*}:=\max \left\{ \widetilde{A}_{K}:2\leq K\leq \nu
+1\right\} $, we have $c_{A,K}<m_{A}$ for every $A>A_{*}$ and $K=2,...,\nu
+1 $, so that all the $\nu $ weak solutions $u_{K}$ with $K=2,...,\nu +1$
are nonradial, since otherwise Lemma \ref{LEM: mr} would yield the
contradiction $c_{A,K}=I\left( u_{K}\right) \geq m_{A}$. This also implies $%
u_{K_{1}}\neq u_{K_{2}}$ for $K_{1}\neq K_{2}$, since $H_{K_{1}}\cap
H_{K_{2}}=H_{\mathrm{r}}$ (see \cite[Lemma 2.1]{R}). \endproof


\end{document}